\newcount\ite\ite=1\def\0{\global\ite=1\1}
\def\1{\item{\rm(\romannumeral\the\ite)}\advance\ite1\quad}

\documentclass[12pt]{amsart}
\usepackage{setspace}

\usepackage{amsmath, amsfonts, latexsym, amssymb, graphicx, multirow, tikz}
\usepackage[all]{xy}
\usepackage[colorlinks, linkcolor=blue, citecolor=magenta, urlcolor=cyan]{hyperref}

\font\teneufm=eufm10 scaled \magstep1
\font\seveneufm=eufm7 scaled \magstep1
\font\fiveeufm=eufm5  scaled \magstep1
\newfam\eufmfam

\textfont\eufmfam=\teneufm
\scriptfont\eufmfam=\seveneufm
\scriptscriptfont\eufmfam=\fiveeufm

\newfam\msbfam
\font\tenmsb=msbm10 scaled \magstep1  \textfont\msbfam=\tenmsb
\font\sevenmsb=msbm7 scaled \magstep1 \scriptfont\msbfam=\sevenmsb
\font\fivemsb=msbm5 scaled \magstep1  \scriptscriptfont\msbfam=\fivemsb

\def\dd#1{\raise1.5pt\hbox{$\,\partial\!$}/\raise-2.5pt\hbox{$\!\partial#1\,$}}

\def\tilde{\widetilde}

\def\5#1{{\mathcal #1}}

\def\CC{{\mathbb C}}

\def\ra{\rightarrow}

\def\GL{\mathop{\rm GL}\nolimits}
\def\SL{\mathop{\rm SL}\nolimits}

\def\Ann{\mathop{\rm Ann}\nolimits}

\def\Hess{\mathop{\rm Hess}\nolimits}
\def\Frac{\mathop{\rm Frac}\nolimits}

\def\Cat{\mathop{\rm Cat}\nolimits}
\def\Soc{\mathop{\rm Soc}\nolimits}
\def\rank{\mathop{\rm rank}\nolimits}

\newcommand\co{\colon \thinspace} 

 \def\HollowBoxx #1#2#3{{\dimen0=#1 \advance\dimen0 by -#2
       \dimen1=#1 \advance\dimen1 by #3
        \vrule height 0pt depth #3 width #2
       \hskip -#3
       \vrule height #1 depth #3 width #3}}
 \def\LeftContraction{\mathord{\kern1.45pt \HollowBoxx{6pt}{3.5pt}{.4pt}}\,}

 \def\HollowBox #1#2#3{{\dimen0=#1 \advance\dimen0 by -#3
       \dimen1=#1 \advance\dimen1 by #3
        \vrule height #1 depth #3 width #3
        \vrule height 0pt depth #3 width #2
        \hskip -#3}}
 \def\RightContraction{\mathord{\, \HollowBox{6pt}{3.1pt}{.4pt}} \kern1.6pt}

\def\qed{{\hfill $\Box$}}
\newtheorem{theorem}{THEOREM}[section]
\newtheorem{corollary}[theorem]{Corollary}
\newtheorem{lemma}[theorem]{Lemma}
\newtheorem{proposition}[theorem]{Proposition}
\newtheorem{conjecture}[theorem]{Conjecture}

\theoremstyle{definition}

\theoremstyle{remark}
\newtheorem{example}[theorem]{Example}
\newtheorem{remark}[theorem]{Remark}

\makeatletter
\def\blfootnote{\xdef\@thefnmark{}\@footnotetext}
\makeatother

\begin{document}

\title[Associated Forms in Classical Invariant Theory]{Associated Forms
\vspace{0.2cm}\\
in Classical Invariant Theory}
\author[Alper]{Jarod Alper}
\author[Isaev]{Alexander Isaev}

\address[Alper]{Mathematical Sciences Institute\\
Australian National University\\
Canberra, ACT 0200, Australia}
\email{jarod.alper@anu.edu.au}

\address[Isaev]{Mathematical Sciences Institute\\
Australian National University\\
Canberra, ACT 0200, Australia}
\email{alexander.isaev@anu.edu.au}

\begin{abstract}
It was conjectured in the recent article {\rm \cite{EI}} that all absolute classical invariants of forms of degree $m\ge 3$ on $\CC^n$ can be extracted, in a canonical way, from those of forms of degree $n(m-2)$ by means of assigning every form with non-vanishing discriminant the so-called associated form. This surprising conjecture was confirmed in {\rm \cite{EI}} for binary forms of degree $m \le 6$ and ternary cubics. In the present paper, we settle the conjecture in full generality.  In addition, we propose a stronger version of this statement and obtain evidence supporting it.
\end{abstract}
\maketitle

\thispagestyle{empty}

\pagestyle{myheadings}
\markboth{Jarod Alper and Alexander Isaev}{Associated Forms in Classical Invariant Theory}\blfootnote{{\bf Mathematics Subject Classification:} 14L24, 13A50, 13H10, 32S25}\blfootnote{{\bf Keywords:} Classical Invariant Theory, Milnor algebras, Artinian Gorenstein algebras, isolated hypersurface singularities.}

\setcounter{section}{0}

\section{Introduction}\label{intro}
\setcounter{equation}{0}

The paper is devoted to a new construction in classical invariant theory proposed in the recent article \cite{EI}. Let ${\mathcal Q}_n^m$ be the vector space of forms of degree $m$ on $\CC^n$ with $n\ge 2$, $m\ge 3$. For a form $Q\in{\mathcal Q}_n^m$ having a non-zero discriminant, consider the Milnor algebra $M(Q)$ of the isolated singularity at the origin of the hypersurface in $\CC^n$ defined by $Q$. As explained in \cite{EI}, the algebra $M(Q)$ gives rise to a form on the quotient space ${\mathfrak m}/{\mathfrak m}^2$ with values in the (one-dimensional) socle $\Soc(M(Q))$ of $M(Q)$, where ${\mathfrak m}$ is the maximal ideal of $M(Q)$. By making a canonical choice of coordinates in ${\mathfrak m}/{\mathfrak m}^2$ and $\Soc(M(Q))$, one obtains a form ${\bf Q}$ of degree $n(m-2)$ on  $\CC^n$.  The form ${\bf Q}$ is called the {\it associated form}\, of $Q$. 

It was shown in \cite{EI} that for certain values of $n$ and $m$ one can recover all absolute classical invariants of forms of degree $m$ on $\CC^n$ (i.e.~all $\GL(n,\CC)$-invariant rational functions on ${\mathcal Q}_n^m$) from those of forms of degree $n(m-2)$ on $\CC^n$ by evaluating the latter on associated forms. Thus, for such $n$ and $m$ the invariant theory of forms in ${\mathcal Q}_n^m$ can be extracted, in a canonical way, from that of forms in ${\mathcal Q}_n^{n(m-2)}$, at least at the level of absolute invariants. This surprising result is known to hold for binary forms of degrees $3\le m\le 6$ as well as for ternary cubics ($n=m=3$), and the proofs rely on explicit lists of generators of the corresponding algebras of invariants. For binary quintics ($n=2$, $m=5$) and binary sextics ($n=2$, $m=6$) the arguments are quite involved and require computer-assisted calculations.

The above facts motivated a conjecture, posed in \cite{EI}, asserting that an analogous statement holds true for all $n$ and $m$. In this paper we show that this conjecture is indeed correct. Unlike the results in \cite{EI}, our approach does not depend on explicit descriptions of the algebras of invariants. Instead, it is based on a systematic study of the map 
$$
\Phi \co X^m_n\ra{\mathcal Q}_n^{n(m-2)},\quad Q\mapsto {\bf Q},
$$
where $X^m_n\subset{\mathcal Q}_n^m$ consists of all forms with non-vanishing discriminant.

The paper is organized as follows. In Section \ref{sect1}, we give a detailed definition of associated form and state the conjecture proposed in \cite{EI} (see Conjecture \ref{conj1}). We also formulate a stronger, more natural version of it in Conjecture \ref{conj2}. As explained in  Section \ref{sect1}, confirming the stronger conjecture would provide a way for extracting a complete system of biholomorphic invariants of homogeneous isolated hypersurface singularities from their Milnor algebras. Producing invariants of this kind is motivated by the so-called {\it reconstruction problem}\, in singularity theory, which is the problem of finding an effective proof of the well-known Mather-Yau theorem (see \cite{MY}). 

In Section \ref{sect2}, we show that two forms are linearly equivalent if and only if their associated forms are linearly equivalent (see Theorem \ref{lineq}).  This statement provides evidence in support of Conjecture \ref{conj2} since the conjecture implies that the associated form ${\bf Q}$ defines the form $Q$ up to linear equivalence. The proof of  Theorem \ref{lineq} is based on an explicit formula, which is of independent interest, that provides a Macaulay inverse system for an arbitrary Artinian Gorenstein algebra (see Proposition \ref{idealdeterms}). The formula implies that the associated form ${\bf Q}$ is in fact an inverse system for the Milnor algebra $M(Q)$; Theorem \ref{lineq} then follows from the Mather-Yau theorem. 

Conjecture \ref{conj1} is confirmed in Section \ref{sect3} (see Theorem \ref{solconj1}). Firstly, we reduce the conjecture to proving that the range $\Phi(X_n^m)$ intersects the subset of stable forms in ${\mathcal Q}_n^{n(m-2)}$ (see Proposition \ref{genappr}) and, secondly, obtain that this intersection is always non-empty (see Proposition \ref{discrimassoc}). Note that in Proposition \ref{discrimassoc} we in fact show  that $\Phi(X_n^m)$ intersects $X_n^{n(m-2)}$, which allows us to introduce second associated forms (see Remark \ref{iterations}). Theorem \ref{solconj1} then implies that for all pairs $n,m$, except $n=2,m=3$, all absolute invariants of forms in ${\mathcal Q}_n^m$ can also be canonically recovered from those of forms in ${\mathcal Q}_n^{n(n(m-2)-2)}$.

Finally, in Section \ref{sect4} we strengthen results of Section \ref{sect3} for the case of binary forms. First of all, we observe that there is a relative $\GL(2,\CC)$-invariant, namely the catalecticant, that does not vanish for all associated forms (see Proposition \ref{catal}). In particular, every element of $\Phi(X_2^m)$ is semi-stable. Furthermore, in Proposition \ref{range} we give a description of all associated forms that are not stable, which implies that the orbit of every point in $\Phi(X_2^m)$ is closed in the affine variety of binary forms of degree $2(m-2)$ with non-zero catalecticant (see Corollary \ref{closerorbit}). These facts yield an easier proof of Theorem \ref{solconj1} for binary forms and lead to a more detailed variant of Conjecture \ref{conj2} in this case (see\linebreak Conjecture \ref{conj3}).

{\bf Acknowledgement.} This work is supported by the Australian Research Council.

\section{The associated form and main conjectures}\label{sect1}
\setcounter{equation}{0}

\subsection{Preliminaries}
Let ${\mathcal Q}_n^m$ be the vector space of forms of degree $m$ on $\CC^n$, where $n\ge 2$. The dimension of this space is given by the well-known formula
\begin{equation}
\dim_{\CC}{\mathcal Q}_n^m=\left(
\begin{array}{c}
m+n-1\\
m
\end{array}
\right).\label{dimmathbfq}
\end{equation}
To every non-zero $Q\in{\mathcal Q}_n^m$ we associate the hypersurface
$$
V_Q:=\{z=(z_1,\dots,z_n)\in\CC^n:Q(z)=0\}
$$
and consider it as a complex space with the structure sheaf induced by $Q$. The singular set of $V_Q$ is then the critical set of $Q$. In particular, if $m\ge 2$ the hypersurface $V_Q$ has a singularity at the origin. In this paper, we focus on the situation when this singularity is isolated, or, equivalently, when $V_Q$ is smooth away from 0. This occurs if and only if $Q$ is {\it non-degenerate}, i.e.~$\Delta(Q)\ne 0$, where $\Delta$ denotes the discriminant (see \cite{GKZ}, Chapter 13). 

For $m\ge 3$ define
$$
X^m_n:=\{Q\in{\mathcal Q}_n^m:\Delta(Q)\ne 0\}.
$$
Fix $Q\in X^m_n$ and consider the {\it Milnor algebra}\, of the singularity\ of $V_Q$, which is the complex local algebra
$$
M(Q):=\CC[[z_1,\dots,z_n]]/J(Q),
$$
where $\CC[[z_1,\dots,z_n]]$ is the algebra of formal power series in $z_1,\dots,z_n$ with complex coefficients and $J(Q)$ the {\it Jacobian ideal}\, in $\CC[[z_1,\dots,z_n]]$, i.e.~the ideal generated by all the first-order partial derivatives\linebreak $Q_j:=\partial Q/\partial z_j$ of $Q$, $j=1,\dots,n$. Observe that $m\ge 3$ implies $\dim_{\CC}M(Q)>n+1$. On the other hand, since the singularity of $V_Q$ is isolated, the algebra $M(Q)$ is {\it Artinian}, i.e. $\dim_{\CC}M(Q)<\infty$. Therefore, $Q_1,\dots,Q_n$ is a system of parameters in $\CC[[z_1,\dots,z_n]]$. Since $\CC[[z_1,\dots,z_n]]$ is a regular local ring and therefore Cohen-Macaulay,  $Q_1,\dots,Q_n$ is a regular sequence in $\CC[[z_1,\dots,z_n]]$. This yields that $M(Q)$ is a complete intersection.

For convenience, from now on we utilize another realization of the Milnor algebra. Namely, we write
$$
M(Q)=\CC[z_1,\dots,z_n]/{J}(Q),
$$
where we abuse notation by also using $J(Q)$ to refer to the ideal in $\CC[z_1,\dots,z_n]$ generated by $Q_j$, $j=1,\dots,n$.

Let ${\mathfrak m}$ denote the maximal ideal of ${M}(Q)$, which consists of all elements represented by polynomials in $\CC[z_1,\dots,z_n]$ vanishing at the origin. Since ${ M}(Q)$ is a complete intersection, \cite{B} implies that ${ M}(Q)$ is a {\it Gorenstein}\, algebra. This means that the socle of ${M}(Q)$, defined as $\Soc({M}(Q)):=\{x\in{\mathfrak m}: x\,{\mathfrak m}=0\}$, is a one-dimensional vector space over $\CC$  (see, e.g.~\cite{Hu}, Theorem 5.3). In fact, $\Soc({M}(Q))$ is spanned by the element $e_0$ represented by the Hessian $\Hess(Q)$ of $Q$ (see, e.g.~\cite{Sa}, Lemma 3.3). Notice that $\Hess(Q)\in{\mathcal Q}_n^{n(m-2)}$. 

The maximal ideal ${\mathfrak m}$ is nilpotent and we denote by $\nu$ the {\it socle degree}\, of ${M}(Q)$, i.e.~the largest among all integers $\eta$ for which ${\mathfrak m}^{\eta}\ne 0$. Clearly, one has $\Soc({M}(Q))={\mathfrak m}^{\nu}$. On the other hand, for every integer $j\ge 0$, let $L_j$ be the subspace of elements of ${M}(Q)$ represented by forms in ${\mathcal Q}_n^j$. These subspaces constitute a grading on ${M}(Q)$, i.e. 
\begin{equation}  
{M}(Q)=\bigoplus_{j=0}^{\infty}L_j,\quad L_iL_j\subset L_{i+j}.\label{grading}
\end{equation}
We then have
$$
{\mathfrak m}=\bigoplus_{j=1}^{\infty}L_j
$$
and $\Soc({M}(Q))=L_d$, with $d:=\max\{j:L_j\ne 0\}$. Notice that the grading $\{L_j\}$ is {\it standard}, that is, $L_j=L_1^j$ for all $j\ge 1$, which implies $d=\nu$. Furthermore, since $\Soc({M}(Q))$ is spanned by an element represented by a form in ${\mathcal Q}_n^{n(m-2)}$, one also has $d=n(m-2)$. We thus see that $\nu=n(m-2)$ and the subspace
$$
{\mathcal Q}_n^{n(m-2)-m+1}Q_1+\dots+{\mathcal Q}_n^{n(m-2)-m+1}Q_n\subset{\mathcal Q}_n^{n(m-2)}\label{subspace}
$$
has codimension 1, with the line spanned by $\Hess(Q)$ being complementary to it.

For future reference, we also note several properties of the {\it Hilbert function}
$$
H(t):=\sum_{j=0}^{n(m-2)}\dim_{\CC}L_j\cdot t^j
$$
of the grading $\{L_j\}$. Namely, recalling that $\dim_{\CC}{\mathcal Q}_n^j$ is given by formula (\ref{dimmathbfq}), we have
\begin{equation}
\begin{array}{ll}
\dim_{\CC} L_j=\dim_{\CC}{\mathcal Q}_n^j,&\hbox{for $j=0,\dots,m-2$,}\\
\vspace{-0.1cm}\\
\dim_{\CC} L_{m-1}=\dim_{\CC}  {\mathcal Q}_n^{m-1}-n,&\\
\vspace{-0.1cm}\\
\dim_{\CC} L_{n(m-2)-j}=\dim_{\CC} L_j,&\hbox{for $0\le j\le n(m-2)/2$,}
\end{array}\label{holbertfunction}
\end{equation}
with the last statement being a consequence, for instance, of Proposition 9 in \cite{W}. 

\subsection{The associated form} Let $\omega \co \Soc({M}(Q))\ra\CC$ be the linear isomorphism satisfying $\omega(e_0)=1$. Define ${\bf Q}\in{\mathcal Q}_n^{n(m-2)}$ by the formula
$$
{\bf Q}(z):=\omega\left((z_1e_1+\dots+z_ne_n)^{n(m-2)}\right),
$$
for $z = (z_1, \ldots, z_n) \in \CC^n$, where $e_j$ is the element of ${M}(Q)$ represented by the coordinate function $z_j$, $j=1,\dots,n$. We call ${\bf Q}$ the {\it associated form}\, of $Q$. Observe that ${\mathbf Q}$ is a coordinate representation of the following $\Soc({M}(Q))$-valued function on the quotient ${\mathfrak m}/{\mathfrak m}^2$:
\begin{equation}
x\mapsto y^{n(m-2)},\label{coordfree}
\end{equation}
where $y$ is any element of ${\mathfrak m}$ representing $x\in{\mathfrak m}/{\mathfrak m}^2$.

To give an expanded expression for ${\bf Q}$, observe that if  $k_1,\dots,k_n$ are non-negative integers such that $\sum_{j=1}^nk_j=n(m-2)$, the product $e_1^{k_1}\dots e_n^{k_n}$ lies in $\Soc({M}(Q))$, and thus we have 
$$
e_1^{k_1}\cdots e_n^{k_n}=\mu_{k_1,\dots,k_n}(Q) e_0
$$
 for some $\mu_{k_1,\dots,k_n}(Q)\in\CC$. In terms of the coefficients $\mu_{k_1,\dots,k_n}(Q)$ the form ${\bf Q}$ is written as 
\begin{equation}
{\bf Q}(z)=\sum_{k_1+\cdots+k_n=n(m-2)}\mu_{k_1,\dots,k_n}(Q)\frac{(n(m-2))!}{k_1!\cdots k_n!}
z_1^{k_1}\cdots z_n^{k_n}.\label{assocformexp}
\end{equation}
Notice that each $\mu_{k_1,\dots,k_n}$ is a regular function on the affine algebraic variety $X_n^m$, hence
$$
\mu_{k_1,\dots,k_n}=\frac{P}{\Delta^p}\Big|_{X_n^m},\label{coeffmu}
$$
where $P\in\CC[{\mathcal Q}_n^m]$ and $p$ is a non-negative integer. 

\subsection{The morphism $\Phi$} \label{S:phi}
We now consider the morphism
\begin{equation}
\Phi_n^m  \co X_n^m\ra{\mathcal Q}_n^{n(m-2)},\quad Q\mapsto {\bf Q}.\label{mapphinm}
\end{equation}
of affine algebraic varieties.  Throughout the paper, we denote this morphism simply by $\Phi \co X_n^m \to {\mathcal Q}_n^{n(m-2)}$ if there is no fear of confusion. As stated in the introduction, our motivation for introducing the map $\Phi$ comes from the recent article \cite{EI}, in which an attempt was made to relate the invariant theory of forms in $X_n^m$ to that of their associated forms -- see Section \ref{csmc} for a detailed discussion. 

There is a natural action of $\GL(n,\CC)$ on ${\mathcal Q}_n^m$, namely
\begin{equation}
(C,Q)\mapsto Q_C,\quad  Q_C(z):=Q\left(z\,(C^{-1})^T\right)\label{formqc}
\end{equation}
for $C\in\GL(n,\CC)$, $Q\in{\mathcal Q}_n^m$ and $z=(z_1,\dots,z_n)\in\CC^n$. Two forms that lie in the same $\GL(n,\CC)$-orbit are called {\it linearly equivalent}. 

It was shown in Theorem 2.1 of \cite{EI} and Proposition 5.7 of \cite{FK} that if $Q$ and $\tilde{Q}$ are two linearly equivalent forms in $\mathcal{Q}_n^m$, then the associated forms ${\bf Q}$ and ${\bf \tilde{Q}}$ are also linearly equivalent.  We will now strengthen this statement by establishing the following fundamental equivariance property of $\Phi$. 
  
\begin{proposition}\label{equivariance}\sl For every $Q\in X_n^m$ and $C\in\GL(n,\CC)$ one has
\begin{equation}
\Phi(Q_C)=(\det C)^2\,\Phi(Q)_{(C^{-1})^T}.\label{formequivar}
\end{equation}

\end{proposition} 

\noindent {\bf Proof:} Let linear forms $z_j^*$ on $\CC^n$ be defined by
$$
\left(\begin{array}{c}
z_1^*\\
\vdots\\
z_n^*
\end{array}
\right)=C^{-1}\left(\begin{array}{c}
z_1\\
\vdots\\
z_n
\end{array}
\right),
$$
and ${\rm e}_j$, ${\rm e}_j^*$ the elements of ${M}(Q_C)$ represented by $z_j$, $z_j^*$, respectively, $j=1,\dots,n$. Further, denote by ${\rm e}_0^*$ the element of $\Soc(M(Q_C))$ represented by $\Hess(Q_C)$. Then for all non-negative integers $k_1,\dots,k_n$ such  that $k_1+\dots+k_n=n(m-2)$ we have
\begin{equation}
{\rm e}_1^{*\,k_1}\dots {\rm e}_n^{*\,k_n}=(\det C)^2\,\mu_{k_1,\dots,k_n}(Q) {\rm e}_0^*.\label{relat888}
\end{equation}

Next, let $\omega^* \co \Soc({M}(Q_C))\ra\CC$ be the linear isomorphism satisfying $\omega^*({\rm e}_0^*)=1$. Formula (\ref{assocformexp}) and identities (\ref{relat888}) then imply
\begin{equation}
\omega^*\left((z_1{\rm e}_1^*+\dots+z_n{\rm e}_n^*)^{n(m-2)}\right)=(\det C)^2\,\Phi(Q)(z)\label{relspec1}
\end{equation} 
for all points $z=(z_1,\dots,z_n)\in\CC^n$. On the other hand, for every $z\in\CC^n$ one has
\begin{equation}
\begin{array}{l}
\omega^*\left((z_1{\rm e}_1^*+\dots+z_n{\rm e}_n^*)^{n(m-2)}\right)=\\
\vspace{-0.3cm}\\
\hspace{2cm}\omega^*\left((w_1{\rm e}_1+\dots+w_n{\rm e}_n)^{n(m-2)}\right)=\Phi(Q_C)(w),
\end{array}\label{relspec2}
\end{equation}
where the point $w=(w_1,\dots,w_n)$ is defined by $w:=z\, C^{-1}$. Identities (\ref{relspec1}) and (\ref{relspec2}) now lead to formula (\ref{formequivar}), as required.\qed

\begin{remark}  We note that Proposition \ref{equivariance} yields the useful fact that the constructible set\, $\Phi(X_n^m)\subset{\mathcal Q}_n^{n(m-2)}$ is $\GL(n,\CC)$-invariant.
\end{remark}

\subsection{Classical invariants and the main conjectures}\label{csmc}

In this section, we introduce the main conjectures and summarize the relevant results of \cite{EI}.

We start by recalling the definitions of relative and absolute classical invariants as well as related concepts (see, e.g.~\cite{GIT}, \cite{Muk}, \cite{O} for details).
A {\it relative invariant}\, (or {\it relative classical invariant}) of forms of degree $m$ on $\CC^n$ is a polynomial $I\co {\mathcal Q}_n^m\ra\CC$ such that for any $Q\in{\mathcal Q}_n^m$ and any $C\in\GL(n,\CC)$ one has $I(Q)=(\det C)^{\ell}I(Q_C)$, where $Q_C$ is the form introduced in (\ref{formqc}) and $\ell$ a non-negative integer called the {\it weight}\, of $I$. It follows that $I$ is in fact homogeneous of degree $\ell n/m$. Finite sums of relative invariants constitute the algebra 
\begin{equation} \label{defn-mathcalA}
{\mathcal A}_n^m := \CC[\mathcal{Q}_n^m]^{\SL(n,\CC)}
\end{equation}
of $\SL(n,\CC)$-invariant regular functions on ${\mathcal Q}_n^m$, called the {\it algebra of invariants}\, (or {\it algebra of classical invariants}) of forms of degree $m$ on $\CC^n$. As shown by Hilbert in \cite{hilbert}, this algebra is finitely generated. Note that ${\mathcal{A}}_n^m$ is graded by the weight of the invariant.

A form $Q_0\in{\mathcal Q}_n^m$ is called {\it semi-stable}, if for some non-constant $I\in{\mathcal A}_n^m$ one has $I(Q_0)\ne 0$. We denote the subset of semi-stable forms in ${\mathcal Q}_n^m$ by $({\mathcal Q}_n^m)^{\rm\small ss}$. Clearly, $({\mathcal Q}_n^m)^{\rm\small ss}$ is open  (here and below all topological statements refer to the Zariski topology unless stated otherwise). Furthermore, a semi-stable form $Q_0\in{\mathcal Q}_n^m$ is called {\it stable}, if\linebreak (i) the isotropy subgroup of $Q_0$ under the $\GL(n,\CC)$-action is finite,\linebreak and (ii) there exists a non-constant $I\in{\mathcal A}_n^m$ such that $I(Q_0)\ne 0$\linebreak with the $\GL(n,\CC)$-orbit of $Q_0$ being closed in the affine variety\linebreak $\{Q\in{\mathcal Q}_n^m: I(Q)\ne 0\}$. It then follows that the orbit of $Q_0$ is closed in $({\mathcal Q}_n^m)^{\rm\small ss}$ (see \cite{GIT}, Amplification 1.11 or \cite{Muk}, Lemma 5.40). The subset of all stable forms is open in $({\mathcal Q}_n^m)^{\rm\small ss}$, and we denote it by $({\mathcal Q}_n^m)^{\rm\small s}$.   Since any form with non-zero discriminant is stable (see \cite{GIT}, Proposition 4.2), one has $X_n^m\subset ({\mathcal Q}_n^m)^{\rm\small s}$.  Thus the $\GL(n, \CC)$-orbit of any form in $X_n^m$ is closed in $X_n^m$.

For any two relative invariants $I$ and $\tilde I$, with $\tilde I\not\equiv 0$, the ratio $I/\tilde I$ yields a rational function on ${\mathcal Q}_n^m$, which is defined, in particular, at the points where $\tilde I$ does not vanish. If $I$ and $\tilde I$ have equal weights, this function does not change under the action of $\GL(n,\CC)$, and one calls $I/\tilde I$ an {\it absolute invariant}\, (or {\it absolute classical invariant}) of forms of degree $m$ on $\CC^n$. Such invariants form the field
\begin{equation*} \label{defn-mathcalI}
\mathcal{I}_n^m := \Frac(\mathcal{A}_n^m)_0
\end{equation*}
of degree 0 elements in the fraction field of $\mathcal{A}_n^m$.  By Proposition 1 of \cite{DC}, the field $\mathcal{I}_n^m$ may be identified with the field of $\GL(n,\CC)$-invariant rational functions on $\mathcal{Q}_n^m$.

Finally, we define
\begin{equation} \label{defn-mathfrakI}
{\mathfrak I}_n^m:=\CC[X_n^m]^{\,\GL(n,\CC)}
\end{equation}
to be the algebra of $\GL(n,\CC)$-invariant regular functions on $X_n^m$. Clear\-ly, every function in ${\mathfrak I}_n^m$ is the restriction to $X_n^m$ of an element of ${\mathcal I}_n^m$ of the form $I/\Delta^p$, where  $I$ is a relative invariant and $p$ a non-negative integer.  In other words, there is an identification
$\mathfrak{I}_n^m \simeq ((\mathcal{A}_n^m)_{\Delta})_0$ with the degree 0 elements in the localization of $\mathcal{A}_n^m$ by the discriminant $\Delta$.

In \cite{EI}, the following conjecture was proposed. 

\begin{conjecture}\label{conj1}\sl For every ${\rm I}\in{\mathfrak I}_n^m$ there exists\, ${\mathbf I}\in{\mathcal I}_n^{n(m-2)}$such that the composition ${\mathbf I}\circ\Phi$ is a regular function on $X_n^m$ and coincides\linebreak with ${\rm I}$.   
\end{conjecture}

We will now restate Conjecture \ref{conj1} in different terms. First, note that, by Proposition \ref{equivariance}, if two forms $Q,\,\tilde Q\in X_n^m$ are linearly equivalent, then their associated forms ${\bf Q},\,  {\bf \tilde{Q}} \in{\mathcal Q}_n^{n(m-2)}$ are linearly equivalent. Therefore, for any ${\bf I}\in{\mathcal I}_n^{n(m-2)}$ such that ${\mathbf I}\circ\Phi$ is defined at least at one point of $X_n^m$, this composition is an invariant rational function on $X_n^m$. Denote by ${\mathcal R}_n^m$ the collection of all invariant rational functions on $X_n^m$ obtained in this way. One has ${\mathcal R}_n^m\subset{\mathcal I}_n^m|_{X_n^m}$ (which follows, for example, from Proposition 1 of \cite{DC}). Since every element of ${\mathcal I}_n^m|_{X_n^m}$ can be represented as a ratio of two functions in ${\mathfrak I}_n^m$ (see \cite{Muk}, Proposition 6.2), Conjecture \ref{conj1} is equivalent to the following identity:
$$
{\mathcal R}_n^m={\mathcal I}_n^m|_{X_n^m}.
$$
Thus, Conjecture \ref{conj1} means that one can extract all absolute invariants of forms of degree $m$ from those of forms of degree $n(m-2)$ by applying the latter to associated forms.

For binary cubics ($n=2$, $m=3$) the conjecture is obvious since all forms in $X_2^3$ are pairwise linearly equivalent. For binary quartics ($n=2$, $m=4$) and ternary cubics ($n=3$, $m=3$) it is easy to verify as well (see \cite{EI} and the earlier article \cite{Ea}). Furthermore, in \cite{EI} Conjecture \ref{conj1} was shown to hold for binary quintics ($n=2$, $m=5$) and binary sextics ($n=2$, $m=6$), in which cases the proofs are much harder and utilize computer algebra. These results rely on explicit lists of generators of the corresponding algebras ${\mathcal A}_n^m$ (see, e.g.~\cite{Sy}), and one can observe from the proofs that in fact in each of the above situations a stronger statement takes place. We formulate it as a conjecture for general $n$, $m$ as follows.

\begin{conjecture}\label{conj2}\sl For every ${\rm I}\in{\mathfrak I}_n^m$ there exists ${\mathbf I}\in{\mathcal I}_n^{n(m-2)}$ defined at all points of the set\, $\Phi(X_n^m)\subset{\mathcal Q}_n^{n(m-2)}$ such that ${\mathbf I}\circ\Phi={\rm I}$ on $X_n^m$.
\end{conjecture}
\noindent A more detailed variant of this conjecture will be given in Section \ref{sect4} for the case of binary forms (see Conjecture \ref{conj3}).

Notice that Conjecture \ref{conj1} is a priori weaker than Conjecture \ref{conj2}. Indeed, it may potentially happen that for some ${\rm I}\in{\mathfrak I}_n^m$ there exist ${\mathbf I}\in{\mathcal I}_n^{n(m-2)}$ such that ${\mathbf I}\circ\Phi$ is regular and coincides with ${\rm I}$ on $X_n^m$, but ${\mathbf I}$ is not defined at every point of $\Phi(X_n^m)$ as a rational function on ${\mathcal Q}_n^{n(m-2)}$. As the following simple example shows, an effect of this kind can certainly occur for morphisms of affine algebraic varieties.

\begin{example}\label{ex1}\rm Let $\varphi \co \CC\ra\CC^2$ be the morphism given by $\varphi(z) = (z,z)$
and consider the rational function ${\bf f}(z_1,z_2):=z_1/z_2$ on $\CC^2$. Clearly, the origin in $\CC^2$ lies in the range of $\varphi$ and is a point of indeterminacy for ${\bf f}$. On the other hand, ${\bf f}\circ\varphi\equiv 1$ is a regular function on $\CC$.
\end{example}

\subsection{Relation to singularity theory} \label{S:motivation}

While the conjectures sta\-ted above are interesting from the invariant-theoretic viewpoint, Conjecture \ref{conj2} also has an important implication for singularity theory, which in fact was the main motivation for article \cite{EI}. Namely, if this conjecture were confirmed, it would lead to a way of extracting a complete system of biholomorphic invariants of homogeneous isolated hypersurface singularities in $\CC^n$ from their Milnor algebras. 

Indeed, for two hypersurface singularities $V_{Q}$ and $V_{\tilde{Q}}$ defined by $Q \in X_n^m$ and $\tilde{Q} \in X_n^{\tilde m}$, a biholomorphic equivalence of the germs of $V_{Q}$ and $V_{\tilde Q}$ induces a linear equivalence of $Q$ and $\tilde{Q}$ (by simply examining the linear part of the biholomorphism).  In particular, $m = \tilde{m}$. We now fix $m$ and consider the collection ${\mathcal C}_n^m$ of all absolute invariants ${\mathbf I}\in{\mathcal I}_n^{n(m-2)}$ provided by Conjecture \ref{conj2}. For every ${\mathbf I}\in{\mathcal C}_n^m$ let $F_{\mathbf I}$ be the function that assigns the singularity of $V_Q$ the value ${\mathbf I}({\mathbf Q})$. We note that this value can be calculated directly from the Milnor algebra $M(Q)$ by utilizing function (\ref{coordfree}) defined in a coordinate-free way instead of ${\mathbf Q}$ and regarding ${\mathbf I}$ as a function on $n(m-2)$-forms on ${\mathfrak m}/{\mathfrak m}^2$. By Proposition \ref{equivariance}, every $F_{\mathbf I}$ is a biholomorphic invariant of the homogeneous isolated hypersurface singularities arising from forms in $X_n^m$. Since every non-degenerate form is stable, it follows that the $\GL(n,\CC)$-orbits in $X_n^m$ are separated by elements of ${\mathfrak I}_n^m$ (see, e.g.~\cite{EI}, Proposition 3.1). Therefore, if two homogenous isolated hypersurface singularities have the same value for every function $F_{\mathbf I}$ with ${\mathbf I}\in{\mathcal C}_n^m$, Conjecture \ref{conj2} implies that the two singularities are biholomorphically equivalent. This shows that $\{F_{\mathbf I}\}_{{\mathbf I}\in{\mathcal C}_n^m}$ is a complete system of biholomorphic invariants of homogeneous isolated hypersurface singularities.  

In \cite{EI}, invariants of this kind were introduced in the more general setting of quasi-homogeneous singularities, with the motivation coming from the well-known {\it reconstruction problem}\, in singularity theory. This problem arises from the Mather-Yau theorem, which states that the Tjurina algebras of two isolated hypersurface singularities in $\CC^n$ are isomorphic if and only if the singularities are biholomorphically equivalent (see \cite{MY}). For the class of quasi-homogeneous (in particular, homogeneous) singularities this theorem is contained in the earlier paper \cite{Sh}, in which case the Tjurina algebras coincide with the Milnor algebras. The proofs of the theorem given in \cite{Sh}, \cite{MY} are not constructive, and the reconstruction problem asks for an effective way of recovering a hypersurface germ from the corresponding algebra. Confirming Conjecture \ref{conj2} would provide a solution to the reconstruction problem in the homogeneous case in the form of the complete system of biholomorphic invariants
$\{F_{\mathbf I}\}_{{\mathbf I}\in{\mathcal C}_n^m}$.
\medskip

In the forthcoming sections we will study the map $\Phi$. As a consequence, Conjecture \ref{conj1} will be established in full generality and substantial evidence in favor of Conjecture \ref{conj2} will be obtained. We will start with results supporting Conjecture \ref{conj2}. 

\section{The fibers of {$\Phi$}}\label{sect2}
\setcounter{equation}{0}

As explained in Section \ref{S:motivation}, Conjecture \ref{conj2} implies that a form\linebreak $Q\in X_n^m$ is determined up to linear equivalence by its associated form ${\bf Q}$. The main result of the present section yields that this consequence of Conjecture \ref{conj2} is indeed correct. In fact, the following stronger statement holds.

\begin{theorem}\label{lineq}\sl Two forms $Q,\tilde Q\in X_n^m$ are linearly equivalent if and only if their associated forms\, ${\mathbf Q}, {\mathbf {\tilde Q}}$ are linearly equivalent. 
\end{theorem}
The necessity implication of Theorem \ref {lineq} follows from Proposition \ref{equivariance}.  Alternatively, as mentioned in Section \ref{S:phi}, this implication is a consequence of Theorem 2.1 of \cite{EI} and Proposition 5.7 of \cite{FK}. Here we give a simultaneous proof of both implications.
  
 \subsection{Macaulay inverse systems} \label{S:inverse}
To prove Theorem \ref{lineq}, we will establish a general result that provides an explicit formula for an inverse system of any Artinian Gorenstein algebra (see Proposition \ref{idealdeterms}).  We will then apply this result to prove that for any non-degenerate form $Q$, the associated form ${\bf Q}$ is an inverse system for the Milnor\linebreak algebra $M(Q)$.

First, we recall the definition of inverse system.  What follows is valid for any field of characteristic zero, and we fix such a field $k$. Let $J$ be an ideal in $k[x_1,\dots,x_n]$ lying in the ideal generated by $x_1,\dots,x_n$. It is well-known that the quotient $A:=k[x_1,\dots,x_n]/J$ is an Artinian Gorenstein algebra if and only if there exists a polynomial $g\in k[x_1,\dots,x_n]$ satisfying $\Ann(g)=J$, where
$$
\Ann(g):=\left\{f\in k[x_1,\dots,x_n]: f\left(\frac{\partial}{\partial x_1},\dots,\frac{\partial}{\partial x_n}\right)(g)=0\right\}\label{annihilatorform}
$$
is the {\it annihilator}\, of $g$ (see, e.g.~\cite{ER} and the references therein). In this case, the degree of $g$ coincides with the socle degree of $A$. The freedom in choosing $g$ with $\Ann(g)=J$ is fully understood, and any such polynomial is called a {\it Macaulay inverse system}, or simply an {\it inverse system}, for the Artinian Gorenstein quotient $A$. Furthermore, if the ideal $J$ is {\it homogeneous}\, (i.e.~generated by forms), one can choose $g$ to be homogeneous (see, e.g.~\cite{Em}, Proposition 7). We note that the classical correspondence $J\leftrightarrow g$ can be also derived from the {\it Matlis duality}\, (see \cite{SV},\linebreak Section 5.4).

Inverse systems solve the isomorphism problem for Artinian Gorenstein quotients. Namely, two such quotients are isomorphic if and only if their inverse systems are equivalent in a certain sense (see \cite{Em}, Proposition 16 and a more explicit formulation in \cite{ER}, Proposition 2.2). In general, the equivalence relation for inverse systems is hard to analyze. However, in the homogeneous case, under the additional assumption that $k$ is algebraically closed, this relation coincides with linear equivalence of forms and therefore leads to the following criterion: two Artinian Gorenstein quotients by homogeneous ideals are isomorphic if and only if their homogeneous inverse systems are linearly equivalent (see \cite{Em}, Proposition 17).

We will now present an explicit way of obtaining an inverse system for any Artinian Gorenstein quotient. To the best of our knowledge, the result that follows is new. Let $A = k[x_1, \ldots, x_n]/J$ be an Artinian Gorenstein quotient as above with socle degree $\eta$. Define ${\rm e}_1,\dots,{\rm e}_n$ to be the elements of $A$ represented by $x_1,\dots,x_n$, respectively. Fix a linear form $\rho \co A\ra k$ whose kernel is complementary to $\Soc(A)$ and set
\begin{equation}
R(x_1,\dots,x_n):=\sum_{j=0}^{\eta}\frac{1}{j!}\rho\left((x_1{\rm e}_1+\dots+x_n{\rm e}_n)^j\right).\label{formr}
\end{equation}

\begin{proposition}\label{idealdeterms}\sl The polynomial $R$ is an inverse system for $A$.
\end{proposition}

\noindent{\bf Proof:} Notice, first of all, that the ideal $J$ consists of all relations among the elements ${\rm e}_1,\dots,{\rm e}_n$, i.e.~polynomials $f\in k[x_1,\dots,x_n]$ with $f({\rm e}_1,\dots,{\rm e}_n)=0$. 

Next, fix any polynomial $f\in k[x_1,\dots,x_n]$
$$
f=\sum_{0\le i_1,\dots,i_n\le N}a_{i_1,\dots,i_n}x_1^{i_1}\dots x_n^{i_n}
$$ 
and calculate
\begin{equation}
\begin{array}{l}
\displaystyle f\left(\frac{\partial}{\partial x_1},\dots,\frac{\partial}{\partial x_n}\right)(R)=\\
\vspace{-0.3cm}\\
\displaystyle\sum_{0\le i_1,\dots,i_n\le N}a_{i_1,\dots,i_n}\sum_{j= i_1+\dots+i_n}^{\eta}\frac{1}{(j-(i_1+\dots+i_n))!}\times\\
\vspace{-0.45cm}\\
\hspace{2.5cm}\rho\Bigl((x_1{\rm e}_1+\dots+x_n{\rm e}_n)^{j-(i_1+\dots+i_n)}{\rm e}_1^{i_1}\dots {\rm e}_n^{i_n}\Bigr)=\\
\vspace{-0.3cm}\\
\displaystyle\sum_{\ell=0}^{\eta}\frac{1}{\ell!}\rho\Bigl((x_1{\rm e}_1+\dots+x_n{\rm e}_n)^\ell\times\\
\vspace{-0.7cm}\\
\hspace{4cm}\displaystyle\sum_{\mbox{\tiny$\begin{array}{l} 0\le i_1,\dots,i_n\le N,\\\vspace{-0.3cm}\\
i_1+\dots+i_n\le\eta-\ell\end{array}$}}
a_{i_1,\dots,i_n}{\rm e}_1^{i_1}\dots {\rm e}_n^{i_n}\Bigr)=\\
\vspace{-0.3cm}\\
\displaystyle\sum_{\ell=0}^{\eta}\frac{1}{\ell!}\rho\Bigl((x_1{\rm e}_1+\dots+x_n{\rm e}_n)^\ell\,f({\rm e}_1,\dots,{\rm e}_n)\Bigr).
\end{array}\label{diff}
\end{equation}
Formula (\ref{diff}) immediately implies $J\subset\Ann(R)$.

Conversely, let $f\in k[x_1,\dots,x_n]$ be an element of $\Ann(R)$. Then (\ref{diff}) yields
\begin{equation}
\sum_{\ell=0}^{\eta}\frac{1}{\ell!}\rho\Bigl((x_1{\rm e}_1+\dots+x_n{\rm e}_n)^\ell\,f({\rm e}_1,\dots,{\rm e}_n)\Bigr)=0.\label{eq2}
\end{equation}
Collecting the terms containing $x_1^{i_1}\dots x_n^{i_n}$ in (\ref{eq2}) we obtain
\begin{equation}
\rho\Bigl({\rm e}_1^{i_1}\dots {\rm e}_n^{i_n}\,f({\rm e}_1,\dots,{\rm e}_n)\Bigr)=0\label{eq3}
\end{equation}
for all indices $i_1,\dots,i_n$. Since ${\rm e}_1,\dots,{\rm e}_n$ generate $A$, identities (\ref{eq3}) yield 
\begin{equation}
\rho\Bigl(a \, f({\rm e}_1,\dots,{\rm e}_n)\Bigr)=0 \qquad \text{for all $a \in A$}.\label{nondeg}
\end{equation}
Further, since the bilinear form $(a,b)\mapsto\rho(ab)$ is non-degenerate on $A$ (see, e.g.~\cite{He}, p.~11), identity (\ref{nondeg}) implies $f({\rm e}_1,\dots,{\rm e}_n)=0$. Therefore $f\in J$, which shows that $J=\Ann(R)$ as required.\qed
\medskip

Suppose now that the ideal $J \subset k[x_1, \ldots, x_n]$ is homogeneous. Then, letting $A_j$ be the subspace of elements represented by forms of degree $j$, we obtain a standard grading on $A$:
$$
A=\bigoplus_{j=0}^{\infty}A_j,\quad A_iA_j\subset A_{i+j},
$$
with $A_j=A_1^j$ for all $j\ge 1$. Hence $\Soc(A)=A_{\eta}$, and we choose $\rho$ satisfying the condition
$$
\ker\rho=\bigoplus_{j=0}^{\eta-1}A_j.
$$
Proposition \ref{idealdeterms} yields the following corollary.

\begin{corollary}\label{cor1}\sl Let the ideal $J$ be homogeneous and $\rho$ chosen as above. Then 
$$
S(x_1,\dots,x_n):=\rho\left((x_1{\rm e}_1+\dots+x_n{\rm e}_n)^{\eta}\right)
$$
is an inverse system for $A$. In particular, for every $Q\in X_n^m$ the associated form ${\bf Q}$ is an inverse system for ${M}(Q)$.
\end{corollary}

\noindent{\bf Proof:} Observe that in this case formula (\ref{formr}) for the inverse system $R$ reduces to that for the polynomial $S$ up to a constant multiple. To obtain the last statement, extend $\omega$ to all of ${M}(Q)$ by setting it to be zero on $\bigoplus_{j=0}^{n(m-2)-1}L_j$ (see (\ref{grading})).\qed

\medskip

Thus, Proposition \ref{idealdeterms} yields a simple proof of the well-known fact that an Artinian Gorenstein quotient by a homogeneous ideal admits a homogeneous inverse system (see, e.g.~\cite{Em}, Proposition 7) and provides an explicit formula for such an inverse system.

Corollary \ref{cor1} implies the following fact, which will be relevant later.

\begin{corollary}\label{cor2}\sl For every $Q\in X_n^m$ the partial derivatives of the associated form ${\bf Q}$ of order $\ell$ are linearly independent elements of ${\mathcal Q}_n^{n(m-2)-\ell}$ if\, $0\le\ell\le m-2$. 
\end{corollary}

\noindent{\bf Proof:} As shown in the first equation of (\ref{holbertfunction}), the ideal ${J}(Q)$\linebreak does not contain forms of degree less than or equal to $m-2$.  Since ${ J}(Q)=\Ann({\mathbf Q})$, the corollary follows.\qed

\subsection{Application to associated forms} We will now apply Corollary \ref{cor1} to obtain the main result of Section \ref{sect2}. 

\medskip

\noindent{\bf Proof of Theorem \ref{lineq}:} Proposition 17 in \cite{Em} and Corollary \ref{cor1} imply that the associated forms ${\mathbf Q}$ and $\tilde{\mathbf Q}$ are linearly equivalent if and only if the algebras ${M}(Q)$ and ${M}(\tilde Q)$ are isomorphic. By the Mather-Yau theorem, ${M}(Q)$ and ${M}(\tilde Q)$ are isomorphic if and only if the germs of the hypersurfaces $V_{Q}$ and $V_{\tilde Q}$ at the origin are biholomorphically equivalent. As explained in Section \ref{S:motivation}, the germs are biholomorphically equivalent if and only if $Q$ and $\tilde Q$ are linearly equivalent. The proof is complete.\qed

\medskip

Notice that the map $\Phi$ is not injective. Indeed, one has, for example, $\Phi(\lambda Q)=\Phi(Q)$ for all $Q\in X_n^m$ if $\lambda^n=1$. On the other hand,  Theorem \ref{lineq} yields the following property of the fibers of $\Phi$.

\begin{corollary}\label{fibres}\sl Every fiber of\, $\Phi$ consists of pairwise linearly equivalent forms. 
\end{corollary}

\noindent This fact implies that the push-forward $\Phi_*{\rm I}$ of any element ${\rm I}\in{\mathfrak I}_n^m$ is a well-defined continuous function on $\Phi(X_n^m)$ (with respect to the Euclidean topology). However, it is not immediately clear whether $\Phi_*{\rm I}$ extends to an absolute classical invariant on ${\mathcal Q}_n^{n(m-2)}$.

\section{Confirmation of Conjecture \ref{conj1} }\label{sect3}
\setcounter{equation}{0}

In this section, we settle Conjecture \ref{conj1}.

\begin{theorem}\label{solconj1}\sl Conjecture {\rm\ref{conj1}} holds for all $n\ge 2$, $m\ge 3$.
\end{theorem}

\noindent Theorem \ref{solconj1} is a direct consequence of the following two propositions.  

\begin{proposition}\label{genappr}\sl Conjecture\, {\rm \ref{conj1}} holds true for a particular pair $n,m$ with $n\ge 2$, $m\ge 3$ provided there exists a form $Q\in X_n^m$ such that its associated form ${\bf Q}$ is stable.
\end{proposition}

\begin{proposition}\label{discrimassoc} \sl For every pair $n,m$ with $n\ge 2$, $m\ge 3$, there exists $Q\in X_n^m$ such that $\Delta({\mathbf Q})\ne 0$.
\end{proposition}

\noindent {\bf Proof of Theorem \ref{solconj1}:} 
Any form with non-zero discriminant is in fact stable (see \cite{GIT}, Proposition 4.2).  Therefore, Proposition \ref{discrimassoc} implies that the hypothesis of Proposition \ref{genappr} holds for every pair\linebreak $n\ge 2$, $m\ge 3$. \qed
\medskip

Propositions \ref{genappr} and \ref{discrimassoc} will be proved in Sections \ref{S:prop1} and \ref{S:prop2}, respectively.

\subsection{Good quotients and the proof of Proposition \ref{genappr}} \label{S:prop1}
Our proof of Proposition \ref{genappr} relies on geometric invariant theory (GIT) quotients of algebraic varieties by reductive algebraic groups. We briefly introduce the relevant properties of such quotients and refer to Chapters 0 and 1 in \cite{GIT} and Chapter 5 in \cite{Muk} for further details.

Let $X$ be a complex algebraic variety and $G$ be a complex reductive algebraic group acting algebraically on $X$. A {\it good quotient} of $X$ by $G$ (if it exists) is an algebraic variety $Z$ such that there is a surjective affine $G$-invariant morphism $\pi \co X\ra Z$ for which $\pi^* \co \mathcal{O}_Z \to \mathcal{O}_X^G$ is an isomorphism, where $\mathcal{O}_Z$ and $\mathcal{O}_X^G$ are the sheafs of regular functions on $Z$ and $G$-invariant regular functions on $X$, respectively. For the purposes of this exposition, we do not list all the properties of good quotients.  Instead, we concentrate only on the following ones, which will be relevant to our arguments below:
\begin{enumerate}
\item[(P1)] $\pi^* \co \CC[Z]\ra\CC[X]^G$ is an isomorphism;
\item[(P2)] for $x,x'\in X$ one has $\pi(x)=\pi(x')$ if and only if $\overline{G \cdot x}\cap\overline{G \cdot x'}\ne\emptyset$ (where $G \cdot x$ is the $G$-orbit of $x$), and every fiber of $\pi$ contains exactly one closed $G$-orbit (the unique orbit of minimal dimension);
\item[(P3)] if $Y$ is an algebraic variety and $\varphi \co X\ra Y$ is a $G$-invariant morphism, then there exists a unique morphism $\tau_{\varphi} \co Z\ra Y$ such that $\varphi=\tau_{\varphi}\circ\pi$;
\item[(P4)] if $A$ is a $G$-invariant closed subset of $X$, then $\pi(A)$ is closed\linebreak in $Z$;
\item[(P5)] if $U\subset X$ is an open subset satisfying $U=\pi^{-1}(\pi(U))$, then properties (P1)--(P4) hold with the triple $(U,\pi(U),\pi|_{U})$ in\linebreak place of the triple $(X,Z,\pi)$, where $\pi(U)\subset Z$ is open by property (P4).
\end{enumerate}
The good quotient of $X$ by $G$ (which is unique up to isomorphism) is denoted by $X/\hspace{-0.1cm}/G$. If every fiber of $\pi$ consists of a single (closed) orbit, the quotient $X/\hspace{-0.1cm}/G$ is said to be {\it geometric}.

For the proof of Proposition \ref{genappr} we will utilize good quotients in two situations. Firstly, a good quotient exists if $X$ is affine. This result goes back to Hilbert for the case $X=\mathcal{Q}_n^m$, $G=\SL(n,\CC)$ (see \cite{hilbert}) and is due to Nagata and and Mumford in general (see \cite{nagata} and \cite{GIT}, Chapter 1, \S 2). In particular, they proved that the algebra of invariants $\CC[X]^G$ is finitely generated. To construct the quotient explicitly, choose generators $f_1,\dots,f_M$ of $\CC[X]^G$ and set
$$
\pi:= (f_1,\dots,f_M): X\ra \CC^M.
$$
Next, consider the ideal
$
K:=\{g\in\CC[z_1,\dots,z_M]: g\circ\pi\equiv 0\}
$
and let
$
Z:=\{z\in\CC^M: g(z) = 0\,\,\hbox{for all}\,\, g\in K\}.
$
Then $Z$ turns out to be a good quotient with the quotient morphism given by $\pi \co X \to Z$, and therefore one can take $Z$ as a realization of $X/\hspace{-0.1cm}/G$. 

Secondly, good quotients are known to exist for the action of $\GL(n, \CC)$ on $\left({\mathcal Q}_n^m\right)^{\rm\small s}$ for any $n\ge 2$, $m\ge 3$ (see, e.g.~\cite{GIT}, Theorem 1.10). In this case, the variety $\left({\mathcal Q}_n^m\right)^{\rm\small s} /\hspace{-0.3cm}/ \GL(n, \CC)$ is quasi-projective. Since every $\GL(n,\CC)$-orbit is closed in $\left({\mathcal Q}_n^m\right)^{\rm\small s}$, the quotient $\left({\mathcal Q}_n^m\right)^{\rm\small s} /\hspace{-0.3cm}/ \GL(n, \CC)$ is geometric by property (P2).  

Moreover, Theorem 1.10 of \cite{GIT} also implies that a good quotient of the semi-stable locus $\left({\mathcal Q}_n^m\right)^{\rm\small{ss}}$ by $\GL(n,\CC)$ exists as well and is the projective variety defined by the graded algebra $\mathcal{A}_n^m$ (see (\ref{defn-mathcalA})).  The quotient of the stable locus $\left({\mathcal Q}_n^m\right)^{\rm\small s} /\hspace{-0.3cm}/ \GL(n, \CC)$ is an open subvariety of $\left({\mathcal Q}_n^m\right)^{\rm\small{ss}} /\hspace{-0.3cm}/ \GL(n, \CC)$.   It follows that the function field of $\left({\mathcal Q}_n^m\right)^{\rm\small{ss}} /\hspace{-0.3cm}/ \GL(n, \CC)$ (and therefore the function field of $\left({\mathcal Q}_n^m\right)^{\rm\small{s}} /\hspace{-0.3cm}/ \GL(n, \CC)$) is canonically identified with the field of absolute invariants $\mathcal{I}_n^m$.

\medskip

\noindent {\bf Proof of Proposition \ref{genappr}:} Consider the affine good quotient    
$$
Z_1:=X_n^m/\hspace{-0.3cm}/\GL(n,\CC),
$$
and let $\pi_1 \co X_n^m\ra Z_1$ be the corresponding $\GL(n,\CC)$-invariant morphism. Every $\GL(n,\CC)$-orbit is closed in $X_n^m$, hence by property (P2) the quotient $Z_1$ is geometric.

Next, consider the stable locus $\bigl({\mathcal Q}_n^{n(m-2)}\bigr)^{\rm\small s}$.  
The hypothesis of Proposition \ref{genappr} implies that $U:=\Phi^{-1}\Bigl(\bigl({\mathcal Q}_n^{n(m-2)}\bigr)^{\rm\small s}\Bigr)$ is a non-empty open subset of $X_n^m$. Further, by Proposition \ref{equivariance}, $U$ is $\GL(n,\CC)$-invariant. Since $Z_1$ is geometric, we then have $U=\pi_1^{-1}(\pi_1(U))$. Therefore, by property (P5) of good quotients, properties (P1)--(P4) hold for the triple $(U,\pi_1(U),\pi_1|_{U})$, where $\pi_1(U)\subset Z_1$ is open.

Next, consider the quasi-projective good quotient
$$
Z_2:=\bigl({\mathcal Q}_n^{n(m-2)}\bigr)^{\rm\small s}/\hspace{-0.3cm}/\GL(n,\CC),
$$
and let $\pi_2 \co \bigl({\mathcal Q}_n^{n(m-2)}\bigr)^{\rm\small s}\ra Z_2$ be the corresponding $\GL(n,\CC)$-invariant morphism. Proposition \ref{equivariance} together with property (P3) yields that there exists a morphism $\phi \co \pi_1(U)\ra Z_2$ such that the diagram
$$\xymatrix{
U \ar[r]^{\Large{\hspace{-0.1cm}\Phi|_U} \quad} \ar[d]^{\pi_{{}_1}|_U} & \bigl({\mathcal Q}_n^{n(m-2)}\bigr)^{\rm\small s} \ar[d]^{\pi_{{}_2}}\\
\pi_1(U)\ar[r]^{\phi} & Z_2
}$$
commutes. Furthermore, Theorem \ref{lineq} implies that $\phi$ is injective.

Next, since the set $\phi(\pi_1(U))$ is constructible, it contains a subset $W$ that is open in the closed irreducible subvariety ${\mathcal R}:=\overline{\phi(\pi_1(U))}$ of $Z_2$. Let ${\mathcal R}_{\rm sng}$ be the singular set of ${\mathcal R}$. Then $W\setminus {\mathcal R}_{\rm sng}$ is non-empty and open in ${\mathcal R}$ as well, and we choose an open subset $O\subset Z_2$ such that $W\setminus {\mathcal R}_{\rm sng}=O\cap {\mathcal R}$. Clearly, $W\setminus {\mathcal R}_{\rm sng}$ is closed in $O$. Next, choose $V\subset O$ to be an affine open subset intersecting $W\setminus {\mathcal R}_{\rm sng}$. Then the set $\tilde{\mathcal R}:=V\cap(W\setminus {\mathcal R}_{\rm sng})$ is closed in $V$ and we have $\tilde{\mathcal R}=V\cap {\mathcal R}$. Let $\tilde U:=\phi^{-1}(V)=\phi^{-1}(\tilde{\mathcal R})$. By construction
$$
\tilde\phi:=\phi|_{\tilde U} \co \tilde U \to \tilde{\mathcal R}\subset V
$$
is a bijective morphism from the irreducible variety $\tilde U$ onto the smooth irreducible variety $\tilde{\mathcal R}$. It now follows from Zariski's Main Theorem that $\tilde\phi$ in fact establishes an isomorphism between $\tilde U$ and $\tilde{\mathcal R}$ (see \cite{TY}, Corollary 17.4.8). Thus, $\tilde \phi \co \tilde{U} \to V$ is a closed immersion of affine varieties.

We will now finalize the proof of the theorem. Fix ${\rm I}\in{\mathfrak I}_n^m$. By property (P1) of good quotients, there is a unique element $f\in\CC[Z_1]$ satisfying $\pi_1^*f={\rm I}$. Since $\tilde \phi \co \tilde{U} \to V$ is a closed immersion of affine varieties, there exists ${\mathbf f} \in \CC[V]$ such that ${\mathbf f}|_{\tilde{\mathcal R}} = \tilde \phi_*\, f|_{\tilde U}$.  Its pull-back $\pi_2^* {\mathbf f}$ is a $\GL(n, \CC)$-invariant regular function on $\pi_2^{-1}(V)$. The continuation ${\mathbf I}$ of $\pi_2^* {\mathbf f}$ to ${\mathcal Q}_n^{n(m-2)}$ is an absolute classical invariant, i.e.~an element\linebreak of ${\mathcal I}_n^{n(m-2)}$.

Notice that by construction the function ${\mathbf I}$ is defined at every point of the open set $\pi_2^{-1}(V)\subset \bigl({\mathcal Q}_n^{n(m-2)}\bigr)^{\rm\small s}$ and, in particular, at every point of $\Phi(\pi_1^{-1}(\tilde U))=\pi_2^{-1}(\tilde{\mathcal R})\subset \pi_2^{-1}(V)$. Furthermore, on the dense open subset $\pi_1^{-1}(\tilde U)\subset X_n^m$ we have ${\mathbf I}\circ\Phi={\rm I}$. It then follows that the rational function ${\mathbf I}\circ\Phi$ is in fact regular on $X_n^m$ and coincides with ${\rm I}$ everywhere as required.\qed

\subsection{Non-degeneracy of associated forms and the proof of Pro\-p\-osition \ref{discrimassoc}} \label{S:prop2}

In our proof of Proposition \ref{discrimassoc}, we will require a characterization of the non-degeneracy of an associated form. We first observe though that Proposition \ref{discrimassoc} can be verified directly in the case of $n=2$, $m=3$. Indeed, every element of $X_2^3$ is linearly equivalent to the binary cubic $z_1z_2(z_1+z_2)$, and an easy computation shows that the associated form of this cubic is proportional to the binary quadratic $z_1^2-z_1z_2+z_2^2$, whose discriminant clearly does not vanish. Therefore, we may exclude the case $n=2$, $m=3$ from our consideration, i.e.~assume that $m\ge 4$ if $n=2$. 

For every $Q\in{\mathcal Q}_n^m$, we introduce the following linear subspace of ${\mathcal Q}_n^{n(m-2)-1}$:
$$
{\mathcal V}(Q):={\mathcal Q}_n^{n(m-2)-m}Q_1+\dots+{\mathcal Q}_n^{n(m-2)-m}Q_n \subset {\mathcal Q}_n^{n(m-2)-1}. 
$$
Further, consider the following cone of powers of linear forms:
$$
{\mathcal C}:=\{(a_1z_1+\dots+a_nz_n)^{n(m-2)-1}: a_1,\dots,a_n\in\CC\}\subset {\mathcal Q}_n^{n(m-2)-1}.
$$

\begin{lemma}\label{critnondegenassoc1}\sl Let $Q\in X_n^m$. Then $\Delta({\bf Q})\ne 0$ if and only if\, ${\mathcal V}(Q)\cap{\mathcal C}=0$.
\end{lemma}

\noindent{\bf Proof:} For a form ${R}\in{\mathcal Q}_n^{\ell}$ with $\ell\ge 2$, the condition $\Delta(R)=0$ means that the partial derivatives ${R}_1,\dots,{ R}_n$ have a common zero, say $z^0$, away from the origin. It then follows that these derivatives vanish on the complex line in $\CC^n$ spanned by $z^0$. Passing to a linearly equivalent form if necessary, we can assume that this line coincides with the\linebreak $z_n$-axis. It is easy to see that a form has all its first-order partial derivatives vanishing on the $z_n$-axis if and only if it can be written as
$$
\sum_{1\le i\le j\le n-1}z_iz_j{ R}^{ij},\quad\hbox{where ${ R}^{ij}\in {\mathcal Q}_n^{\ell-2}$},
$$
i.e.~if it does not involve the monomials $z_n^{\ell}$ and $z_iz_n^{\ell-1}$, with\linebreak $i=1,\dots,n-1$.

Now let $Q$ be an element of $X_n^m$. Then, by the above argument, $\Delta({\mathbf Q})=0$ implies that ${\mathbf Q}$ is linearly equivalent to a form ${\mathbf{\tilde Q}}$ not involving the monomials $z_n^{n(m-2)}$ and $z_iz_n^{n(m-2)-1}$, with $i=1,\dots,n-1$. By Proposition \ref{equivariance}, the form ${\mathbf{\tilde Q}}$ is the associated form of some $\tilde Q\in X_n^m$ linearly equivalent to $Q$. Let $\tilde e_j$ be the element of ${M}(\tilde Q)$ represented by the coordinate function $z_j$, $j=1,\dots,n$, and $\{\tilde L_j\}$ the standard grading on ${M}(\tilde Q)$ constructed as in (\ref{grading}). Then 
\begin{equation}
\begin{array}{ll}
\tilde e_n^{\,\,n(m-2)}=0, &\\
\vspace{-0.1cm}\\
\tilde e_i\tilde e_n^{\,\,n(m-2)-1}=0, & i=1,\dots, n-1.
\end{array}\label{tildeen}
\end{equation}
Identities (\ref{tildeen}) imply that $\tilde e_n^{\,\,n(m-2)-1}\in\Soc({M}(\tilde Q))=\tilde L_{n(m-2)}$. On the other hand, we have $\tilde e_n^{\,\,n(m-2)-1}\in\tilde L_{n(m-2)-1}$, hence $\tilde e_n^{\,\,n(m-2)-1}=0$, which means $z_n^{\,\,n(m-2)-1}\in{\mathcal V}(\tilde Q)$. Since $Q$ is linearly equivalent to\linebreak $\tilde Q$, it follows that there exists a linear form ${\mathfrak L}\in{\mathcal Q}_n^1$ such that\linebreak ${\mathfrak L}^{n(m-2)-1}\in{\mathcal V}(Q)$, which yields ${\mathcal V}(Q)\cap{\mathcal C}\ne 0$. We have thus shown that $\Delta({\mathbf Q})=0$ implies ${\mathcal V}(Q)\cap{\mathcal C}\ne 0$. 

The opposite implication is obtained by observing that the above argument is in fact reversible. This concludes the proof of the\linebreak lemma.\qed

\medskip

Notice that if $Q \in X^m_n$, then $\mathcal{V}(Q)$ is necessarily a codimension $n$ subspace of $\mathcal{Q}_n^{n(m-2)-1}$. Indeed, by properties (\ref{holbertfunction}) of the Hilbert function of the grading $\{L_j\}$ on ${M}(Q)$ (see (\ref{grading})), one has
$$
\dim_{\CC}L_{n(m-2)-1}=\dim_{\CC}L_1=n,
$$
which means
\begin{equation}
\dim_{\CC}{\mathcal V}(Q)=\dim_{\CC}{\mathcal Q}_n^{n(m-2)-1}-n.\label{formdim}
\end{equation}

In the following lemma, we show that if one can construct a\linebreak (possibly degenerate) form ${\mathsf Q} \in \mathcal{Q}^m_n$ such that  $\mathcal{V}({\mathsf Q}) \cap \mathcal{C} = 0$ and $\mathcal{V}({\mathsf Q}) \subset \mathcal{Q}_n^{n(m-2)-1}$ is of codimension $n$, then there is a deformation $Q$ of ${\mathsf Q}$ with $\Delta(Q) \neq 0$ and $\mathcal{V}(Q) \cap \mathcal{C} = 0$.  Then Lemma \ref{critnondegenassoc1} implies that $\Delta({\bf Q}) \neq 0$, as asserted in Proposition \ref{discrimassoc}.

\begin{lemma}\label{critnondegenassoc}\sl Assume that for some ${\mathsf Q} \in{\mathcal Q}_n^m$ one has
$$
\begin{array}{l}
\hbox{{\rm (i)} $\dim_{\CC}{\mathcal V}({\mathsf Q})=\dim_{\CC}{\mathcal Q}_n^{n(m-2)-1}-n$,}\\
\vspace{-0.1cm}\\
\hbox{{\rm (ii)} ${\mathcal V}({\mathsf Q})\cap{\mathcal C}=0$.}
\end{array}
$$
Then there exists $Q\in X_n^m$ such that $\Delta({\mathbf Q})\ne 0$. 
\end{lemma}

\noindent {\bf Proof:} 
Let $Q^k$ be a sequence in $X_n^m$ converging to ${\mathsf Q}$ in the Euclidean topology of ${\mathcal Q}_n^m$ as $k\ra\infty$. Since ${\mathcal V}({\mathsf Q})\cap{\mathcal C}=0$ by assumption and $\dim_{\CC}{\mathcal V}(Q^k)=\dim_{\CC}{\mathcal V}({\mathsf Q})$ for all $k$ by formula (\ref{formdim}), it follows that ${\mathcal V}(Q^k)\cap{\mathcal C}=0$ for all sufficiently large $k$. Hence, by Lemma \ref{critnondegenassoc1}, the associated form of every such $Q^k$ has non-vanishing discriminant, and the lemma is established.\qed
\medskip

\noindent{\bf Proof of Proposition \ref{discrimassoc}:} 
For every pair $n,m$,  we need to produce a form ${\mathsf Q}$ satisfying the assumptions of Lemma \ref{critnondegenassoc}. In fact, we will define a form ${\mathsf Q}$ for which ${\mathcal V}({\mathsf Q})={\mathcal W}$, where ${\mathcal W}$ is the subspace of ${\mathcal Q}_n^{n(m-2)-1}$ spanned by all monomials in $z_1,\dots,z_n$ of degree
 $n(m-2)-1$ other than $z_1^{n(m-2)-1},\dots,z_n^{n(m-2)-1}$. Clearly, ${\mathcal W}$ is a codimension $n$ subspace in ${\mathcal Q}_n^{n(m-2)-1}$ and ${\mathcal W}\cap{\mathcal C}=0$.

Set
$$
{\mathsf Q}(z):=\left\{
\begin{array}{ll}
\displaystyle \sum_{1\le i<j<k\le n}z_iz_jz_k & \hbox{if $m=3$,}\\
\vspace{-0.1cm}\\
\displaystyle \sum_{1\le i<j\le n}(z_i^{m-2}z_j^2+z_i^2z_j^{m-2}) & \hbox{if $m\ge 4$.}
\end{array}
\right.
$$ 
One immediately observes that ${\mathcal V}({\mathsf Q})\subset{\mathcal W}$. To see that ${\mathcal V}({\mathsf Q})={\mathcal W}$, let ${\tt m}_1,\dots,{\tt m}_K$ be the monomial basis in ${\mathcal W}$, where the monomials are listed with respect to some ordering, with $K:=\dim_{\CC}{\mathcal Q}_n^{n(m-2)-1}-n$. For any elements $P^1,\dots,P^n\in{\mathcal Q}_n^{n(m-2)-m}$ we write the expression $P^1{\mathsf Q}_1+\dots+P^n{\mathsf Q}_n$ as
$$
A^1{\tt m}_1+\dots+A^K{\tt m}_K,
$$
where $A^i$ are linear functions of the coefficients $p^1,\dots,p^N$ of the\linebreak forms $P^1,\dots,P^n$ listed with respect to some ordering, with\linebreak $N:=n\dim_{\CC}{\mathcal Q}_n^{n(m-2)-m}$. Further, we write
$$
\left(
\begin{array}{c}
A^1\\
\vdots\\
A^K
\end{array}
\right)=A
\left(
\begin{array}{c}
p^1\\
\vdots\\
p^N
\end{array}
\right),
$$
where $A$ is a constant $K\times N$-matrix independent of $P^1,\dots,P^n$. We also record the inequality $K\le N$, which follows from identity (\ref{formdim}).\linebreak  In terms of the matrix $A$ the condition ${\mathcal V}({\mathsf Q})={\mathcal W}$ means that\linebreak $\rank(A)=K$. A straightforward albeit tedious calculation, which we omit because of its length, now yields that $A$ is indeed of maximal rank thus proving the proposition.\qed

\begin{remark}\label{iterations} \rm Proposition \ref{discrimassoc} implies that $U:=\Phi^{-1}\bigl(X_n^{n(m-2)}\bigr)$\linebreak is an open dense subset of $X_n^m$. Observe also that one always has $n(m-2)\ge 3$ with the exception of the pair $n=2,m=3$. Hence, excluding this special case, for every $Q\in U$ one can introduce the
 {\it second associated form}\, as $(\Phi_n^{n(m-2)}\circ \Phi_n^m) (Q)=\Phi_n^{n(m-2)}({\mathbf Q})$ (see (\ref{mapphinm})). The arguments of this section then show that Conjecture \ref{conj1} holds for second associated forms in place of associated forms, that is, all absolute invariants of elements in ${\mathcal Q}_n^m$ can be recovered from those of elements in ${\mathcal Q}_n^{n(n(m-2)-2)}$ by evaluating the latter for second associated forms. It would be interesting to investigate whether, considering compositions of $k$ maps $\Phi_n^M$ for suitable values of $M$, one can define a $k$th associated form, with any $k\ge 3$, for every $Q$ in a certain dense open subset of $X_n^m$ and obtain a statement analogous to Conjecture \ref{conj1} for such higher associated forms. For example, to introduce a third associated form, one would need to show that 
$$
(\Phi_n^{n(m-2)}\circ \Phi_n^m)(U)\cap X_n^{n(n(m-2)-2)}\ne\emptyset.
$$  
Notice that $n(m-2)>m$ in all situations except $n=2$, $m=4$\linebreak and $n=3$, $m=3$, while in the latter two cases we have\linebreak $n(m-2)=m$. Thus, with the exclusion of these two special cases, the degree of the $k$th associated form of an element of ${\mathcal Q}_n^m$ (provided such a form could be introduced) would tend to infinity as $k\ra\infty$. It would then follow that in all situations, except the three cases $n=2,m=3$; $n=2,m=4$; $n=3,m=3$, there exists a strictly increasing sequence $\{L_k\}$ of natural numbers such that all absolute invariants of forms in ${\mathcal Q}_n^m$ can be canonically retrieved from those of forms in ${\mathcal Q}_n^{L_k}$ for\linebreak every $k$.
\end{remark}
 
\section{The case of binary forms}\label{sect4}
\setcounter{equation}{0}

In this section, we obtain more specialized results in the case of binary forms ($n=2$). These results yield a simpler proof of Theorem \ref{solconj1} as well as motivate a more precise variant of Conjecture \ref{conj2} for binary forms (see Conjecture \ref{conj3}).

Every non-zero binary form can be represented as a product of linear factors; this property will be extensively used in our arguments below. Note that in terms of linear factors, for a form $Q\in{\mathcal Q}_2^m$ the condition $\Delta(Q)\ne 0$ means that each of the $m$ factors of $Q$ has multiplicity one (see, e.g.~\cite{O}, Theorem 2.39), or, equivalently, that the partial derivatives $Q_1$, $Q_2$ of $Q$ do not have common factors. Furthermore, the semi-stability of $Q$ is equivalent to each factor having multiplicity\linebreak not exceeding $m/2$, and the stability of $Q$ is equivalent to each factor having multiplicity less than $m/2$ (see, e.g. \cite{GIT}, Proposition 4.1 or \cite{Muk}, Propo\-sition 7.9).

\subsection{(Semi-)stability of the associated form}
Let
$$
Q(z)=\sum_{j=0}^{2N}\left(\begin{array}{c}
2N\\
j
\end{array}
\right)a_jz_1^{2N-j}z_2^j
$$
be a binary form of even degree $2N$. The {\it catalecticant}\, of $Q$ is defined as follows:
$$
\Cat(Q):=\det\left(\begin{array}{cccc}
a_0 & a_1 & \dots & a_N\\
a_1 & a_2 & \dots & a_{N+1}\\
\vdots & \vdots & \ddots & \vdots\\
a_N & a_{N+1} & \dots & a_{2N}
\end{array}
\right).
$$
The catalecticant is a relative invariant of degree $N+1$ of forms in ${\mathcal Q}_2^{2N}$ (see \cite{Ell}, \S 209). It is well-known (and in fact easy to prove) that $\Cat(Q)=0$ if and only if the partial derivatives of $Q$ of order $N$ are linearly dependent in ${\mathcal Q}_2^N$ (see, e.g.~\cite{K}, Lemma 6.2).  Moreover, the set where the catalecticant vanishes is the closure of the locus of forms in $\mathcal{Q}_2^{2N}$ expressible as the sum of the $(2N)$th powers of $N$ linear forms (see, e.g. \cite{Ell}, \S 208 or \cite{grace-young}, \S 187).

Notice that for $Q\in X_2^m$ one has ${\mathbf Q}\in{\mathcal Q}_2^{2(m-2)}$, so $\Cat({\mathbf Q})$ is well-defined. Corollary \ref{cor2} for $\ell=m-2$ now implies the following proposition.

\begin{proposition}\label{catal}\sl For every $Q\in X_2^m$ one has $\Cat({\mathbf Q})\ne 0$. 
\end{proposition}

\begin{remark}  Alternatively, this proposition can be proven directly, without relying on the material on inverse systems developed in Section \ref{S:inverse} and, in particular, without appealing to Corollary \ref{cor2}.
\end{remark}
For any $N\ge 1$, define
$$
Y_2^{2N}:=\{Q\in{\mathcal Q}_2^{2N}: \Cat(Q)\ne 0\}.
$$
By Proposition \ref{catal}, the map $\Phi$ is a morphism of affine algebraic varieties
$$
\Phi \co X_2^m\ra Y_2^{2(m-2)}.
$$
Since $Y_2^{2(m-2)}\subset \bigl({\mathcal Q}_2^{2(m-2)}\bigr)^{\rm\small ss}$, 
the associated form ${\bf Q}$ of every form $Q \in X_2^m$ is semi-stable.  However, not every associated form is stable.  Indeed, a quick computation shows that for $Q = z_1^m+z_2^m$ the associated form is
$$
{\bf Q} = \frac{1}{m^2(m-1)^2}{2m-4 \choose m-2} z_1^{m-2} z_2^{m-2},
$$
which is clearly not stable. We will now prove that this example in fact describes the only possibility, up to linear equivalence, for an associated form not to be stable.

\begin{proposition}\label{range}\sl If the associated form\, ${\mathbf Q}$ of a form $Q \in X_2^m$ is not stable, then $Q$ is linearly equivalent to $z_1^m+z_2^m$.  In particular, ${\mathbf Q}$ is linearly equivalent to $z_1^{m-2}z_2^{m-2}$.
\end{proposition}

\noindent{\bf Proof:} If the form ${\mathbf Q}$ is not stable, by Proposition \ref{equivariance} we can assume, without loss of generality, that ${\mathbf Q}$ is divisible by $z_1^{m-2}$, that is, $\mu_{0,2m-4}(Q)=\dots=\mu_{m-3,m-1}(Q)=0$ (see (\ref{assocformexp})). 

First, we claim
\begin{equation}
e_2^{m-1}=0.\label{e2power}
\end{equation}
If $m=3$, there is nothing to prove since in this case $\mu_{0,2}(Q)=0$, so we may assume $m>3$. 
From the conditions $\mu_{0,2m-4}(Q)=\mu_{1,2m-5}(Q)=0$, we have
\begin{equation}
\begin{array}{l}
\displaystyle z_2^{2m-4}=\sum_{j=0}^{m-3} z_1^{m-j-3} z_2^j (\alpha_j Q_1+\beta_j Q_2),\\
\vspace{-0.1cm}\\  
\displaystyle z_1z_2^{2m-5}=\sum_{j=0}^{m-3} z_1^{m-j-3} z_2^j (\gamma_j Q_1+\delta_j Q_2)
\end{array}\label{rel11}
\end{equation}
for some $\alpha_j,\beta_j,\gamma_j,\delta_j\in\CC$. Multiplying the first identity in (\ref{rel11}) by $z_1$, the second one by $z_2$, and comparing the results, we see
\begin{equation}
P^1Q_1+P^2Q_2=0,\label{p1p1}
\end{equation}
where
$$
\begin{array}{l}
\displaystyle P^1:=\sum_{j=0}^{m-3}(\alpha_j z_1^{m-j-2} z_2^j -\gamma_j z_1^{m-j-3} z_2^{j+1}),\\
\vspace{-0.1cm}\\
\displaystyle P^2:=\sum_{j=0}^{m-3}(\beta_j z_1^{m-j-2} z_2^j-\delta_j z_1^{m-j-3} z_2^{j+1}).
\end{array}
$$

Since $Q_1, Q_2\in{\mathcal Q}_2^{m-1}$ do not have common linear factors and $P^1$, $P^2$ lie in ${\mathcal Q}_2^{m-2}$, identity (\ref{p1p1}) implies $P^1=P^2=0$. In particular, we obtain $\alpha_0=\beta_0=0$. Division of the first equation in (\ref{rel11}) by $z_2$ then implies
\begin{equation}
z_2^{2m-5}=\sum_{j=0}^{m-4} z_1^{m-j-4}z_2^j (\alpha_{j+1} Q_1+\beta_{j+1} Q_2),\label{idents7788}
\end{equation}
hence
\begin{equation}
e_2^{2m-5}=0.\label{e2power1}
\end{equation}

If $m=4$, identity (\ref{e2power1}) coincides with (\ref{e2power}), otherwise the condition $\mu_{2,2m-6}(Q)=0$ leads to
\begin{equation}
\displaystyle z_1^2z_2^{2m-6}=\sum_{j=0}^{m-3} z_1^{m-j-3}z_2^j (\gamma_j' Q_1+\delta_j' Q_2)
\label{rel22}
\end{equation}
for some $\gamma_j',\delta_j'\in\CC$. Multiplying identity (\ref{idents7788}) by $z_1^2$, identity (\ref{rel22}) by $z_2$, and comparing the polynomial coefficients at $Q_1$, $Q_2$ as above, we obtain $\alpha_1=\beta_1=0$. Dividing (\ref{idents7788}) by $z_2$ then implies $e_2^{2m-6}=0$. Repeating this argument sufficiently many times yields (\ref{e2power}) as claimed.

Next, write
\begin{equation}
\begin{array}{l}
Q_1=az_2^{m-1}+R,\\
\vspace{-0.1cm}\\
Q_2=bz_2^{m-1}+S,
\end{array}\label{eq1}
\end{equation}
where $a,b\in\CC$ and $R,S$ are forms in ${\mathcal Q}_2^{m-1}$ not involving the monomial $z_2^{m-1}$. Since $e_2^{m-1}=0$, there exist $\alpha,\beta\in\CC$ with $|\alpha|+|\beta|>0$ such that $\alpha R+\beta S=0$.  If $\alpha=0$ then $S= 0$, hence $b\ne 0$ and $Q$ is linearly equivalent to $z_1^m+z_2^m$. If $\beta=0$ then $R= 0$, which is impossible since in this case $Q$ has a linear factor of multiplicity greater than one. Thus, we can assume that $\alpha\ne0$, $\beta\ne 0$, and therefore $S=\gamma R$ for some $\gamma\ne 0$.

Let
$$
R=a_1z_1z_2^{m-2}+\dots+a_{m-1}z_1^{m-1},
$$ 
with $a_1,\dots,a_{m-1}\in\CC$. Then, integrating each of equations (\ref{eq1}) and comparing the results, we see
\begin{equation}
a_j=\left(\begin{array}{c} m-1\\ j\end{array}\right)\frac{a}{\gamma^j},\quad j=1,\dots,m-1.\label{eqq2}
\end{equation}
Now, the first identity in (\ref{eq1}) and relations (\ref{eqq2}) yield
$$
Q_1=a\left(\frac{z_1}{\gamma}+z_2\right)^{m-1}.
$$
Therefore, we have $a\ne 0$ and
$$
Q=\frac{a\gamma}{m}\left(\frac{z_1}{\gamma}+z_2\right)^m+c z_2^m,
$$
where $c:=(b-a\gamma)/m$. Hence $c\ne 0$ and $Q$ is linearly equivalent to $z_1^m+z_2^m$ as required.\qed
\medskip

Proposition \ref{range} yields the following fact.

\begin{corollary}\label{closerorbit} \sl For every $Q\in X_2^m$ the $\GL(2,\CC)$-orbit of\, ${\bf Q}$ is closed in $Y_2^{2(m-2)}$.
\end{corollary}

\noindent {\bf Proof:} Firstly, since the orbit of any stable form is closed in the open set of semi-stable forms, it follows that if ${\mathbf Q}$ is stable, its orbit is closed in $\bigl({\mathcal Q}_2^{2(m-2)}\bigr)^{\rm\small ss}$, hence in $Y_2^{2(m-2)}$. Secondly, it is straightforward to observe that the orbit of $z_1^{m-2}z_2^{m-2}$ is closed in $\bigl({\mathcal Q}_2^{2(m-2)}\bigr)^{\rm\small ss}$ as well.\qed

\subsection{Alternative proof of Theorem \ref{solconj1} for binary forms and a refined conjecture}

Corollary \ref{closerorbit} leads to a more explicit proof of Theorem \ref{solconj1} for the case $n=2$, which avoids the use of Proposition \ref{discrimassoc}. Although the proof still relies on the sufficiency implication of Theorem \ref{lineq}, we note that for binary forms this implication can be obtained directly, without utilizing inverse systems.
\medskip

\noindent{\bf Alternative Proof of Theorem \ref{solconj1} for $n=2$:} Consider the affine good quotients    
$$
\begin{array}{l}
Z_1:=X_2^m/\hspace{-0.3cm}/\GL(2,\CC),\\
\vspace{-0.1cm}\\
Z_2:=Y_2^{2(m-2)}/\hspace{-0.3cm}/\GL(2,\CC),
\end{array}
$$
and let $\pi_1 \co X_2^m\ra Z_1$, $\pi_2 \co Y_2^{2(m-2)}\ra Z_2$  be the corresponding $\GL(2,\CC)$-invariant morphisms. As in the proof of Proposition \ref{genappr}, there exists a morphism $\phi \co Z_1\ra Z_2$ such that the diagram
$$\xymatrix{
X_2^m \ar[r]^{\hspace{-0.3cm}\Phi} \ar[d]^{\pi_{{}_1}} &Y_2^{2(m-2)} \ar[d]^{\pi_{{}_2}}\\
Z_1	\ar[r]^{\hspace{-0.3cm}\phi} & Z_2
}$$
commutes. Now, Corollary \ref{closerorbit} and Theorem \ref{lineq} together with property (P2) of good quotients imply that the map $\phi$ is injective. The rest of the argument proceeds as in the proof of Proposition \ref{genappr}.
\qed

\begin{remark}\label{remdiscr} \rm While the above argument does not utilize Proposition \ref{discrimassoc}, this proposition is still of independent interest (see Remark \ref{iterations}). In the case of binary forms, one can give a more streamlined proof of this statement compared to the one presented in Section \ref{S:prop2}. Indeed, one can observe directly from Lemma \ref{critnondegenassoc1} that the discriminant of the associated form of
$$
Q_0(z):=\left\{
\begin{array}{ll}
\displaystyle z_1^4+z_1^2z_2^2+z_2^4 & \hbox{if $m=4$,}\\
\vspace{-0.1cm}\\
\displaystyle z_1^{m-2}z_2+z_1z_2^{m-1} & \hbox{if $m\ne 4$}
\end{array}
\right.
$$ 
does not vanish. In particular, Lemma \ref{critnondegenassoc} is not required for obtaining Proposition \ref{discrimassoc} in the case $n=2$.  
\end{remark}

The proof of Theorem \ref{solconj1} given in this section motivates a more detailed variant of Conjecture \ref{conj2} in the case of binary forms. In what follows, ${\mathfrak I}_2^m:=\CC[X_2^m]^{\GL(2,\CC)}$ is the algebra of $\GL(2,\CC)$-invariant regular functions on $X_2^m$ (see (\ref{defn-mathfrakI})). 
  
\begin{conjecture}\label{conj3}\sl For every ${\rm I}\in{\mathfrak I}_2^m$ there exists an invariant regular function ${\mathbf I}$ on $Y_2^{2(m-2)}$ satisfying\, ${\mathbf I}\circ\Phi={\rm I}$ on $X_2^m$.
\end{conjecture}

We claim that Conjecture \ref{conj3} is equivalent to the statement that the morphism $\phi \co Z_1 \to Z_2$ in the above proof is a closed immersion. Indeed, the conjecture is the statement that every function ${\rm I} \in  {\mathfrak I}_2^m$ extends, under the morphism $\Phi$, to a function in $\CC[Y_2^{2(m-2)}]^{\GL(2,\CC)}$. On the other hand, $\phi \co Z_1 \to Z_2$ is a closed immersion if and only if $\phi^* \co \CC[Z_2] \to \CC[Z_1]$ is surjective. The claim now follows from the identifications
$$
\begin{aligned}
\CC[Z_1] &\simeq {\mathfrak I}_2^m, \\
\CC[Z_2] &\simeq \CC[Y_2^{2(m-2)}]^{\GL(2,\CC)},
\end{aligned}
$$
which are a consequence of property (P1) of good quotients.

Recall that the alternative proof of Theorem \ref{solconj1} given above establishes the injectivity of the morphism $\phi \co Z_1 \to Z_2$.  Therefore, the condition for $\phi$ to be a closed immersion would follow from Zariski's Main Theorem provided one could show that $\phi(Z_1)$ is a closed normal subvariety of $Z_2$.

Further, we observe that there are identifications 
$$ 
\begin{aligned}
{\mathfrak I}_2^m &\simeq ((\mathcal{A}_2^m)_{\Delta})_0, \\
 \CC[Y_2^{2(m-2)}]^{\GL(2,\CC)} & \simeq  ((\mathcal{A}_2^{2(m-2)})_{\Cat})_0,
 \end{aligned}
$$
where we recall from (\ref{defn-mathcalA}) that $\mathcal{A}_2^m :=\CC[\mathcal{Q}_2^m]^{\SL(2,\CC)}$.  Therefore, Conjecture \ref{conj3} can be also formulated by requiring that every absolute classical invariant of the form
\begin{equation}
{\rm I}=\frac{I}{\Delta^p}\label{formdelta}
\end{equation}
on $\mathcal{Q}_2^m$ for some $I\in\CC[{\mathcal Q}_2^m]$ and $p \ge 0$ extends, under the morphism $\Phi$, to an absolute classical invariant of the form
\begin{equation}
{\mathbf I}=\frac{{\mathtt I}}{\Cat^q}\label{formcat}
\end{equation}
on $\mathcal{Q}_2^{2(m-2)}$ for some ${\mathtt I}\in\CC[{\mathcal Q}_2^{2(m-2)}]$ and $q \ge 0$.

Finally, we remark that results of article \cite{EI} confirm Conjecture \ref{conj3} for binary forms of degrees $3\le m\le 6$ (see also \cite{Ea}).  Indeed, the computations performed in \cite{EI} explicitly verify that every absolute classical invariant ${\rm I}$ as in (\ref{formdelta}) extends under $\Phi$ to an absolute classical invariant ${\mathbf I}$ as in (\ref{formcat}).

\end{document}